\documentclass[12pt,letter]{amsart}

\usepackage{amsfonts}
\usepackage{amssymb}
\newcommand\nc{\newcommand}

\newtheorem{theorem}{Theorem}[section]
\newtheorem{prop}[theorem]{Proposition}
\newtheorem{prblm}[theorem]{Problem}
\newtheorem{defin}[theorem]{Definition}
\newtheorem{caution}[theorem]{Caution}
\newtheorem{jtheorem}[theorem]{Theorem}

\newtheorem{remark}[theorem]{Remark}
\newtheorem{restroom}[theorem]{Restroom}
\newtheorem{lemma}[theorem]{Lemma}
\newtheorem{construction}[theorem]{Construction}
\newtheorem{corollary}[theorem]{Corollary}
\newtheorem{example}[theorem]{Example}
\newtheorem{conclusion}[theorem]{Conclusion}
\newtheorem{triviality}[theorem]{Triviality}
\newtheorem{proto}[theorem]{Prototype Quasifibration}
\newtheorem{cauex}[theorem]{Cautionary Example}
\newtheorem{propositiondef}[theorem]{Proposition-Definition}
\newtheorem{subth}{Nuisance}[theorem]
\newtheorem{ssubth}{ }[subth]
\newtheorem{conjecture}[theorem]{Conjecture}

\newcommand{\field}[1]{\ensuremath{\mathbb{#1}}}
\newcommand{\N}{\field{N}}

\nc\tri[1]{\begin{triviality}
\label{#1}}
\nc\rest[1]{\begin{restroom}
\label{#1}}
\nc\conj[1]{\begin{conjecture}
\label{#1}}
\nc\prodef[1]{\begin{propositiondef}
\label{#1}}
\nc\prt[1]{\begin{proto}
\label{#1}}
\nc\lem[1]{\begin{lemma}
\label{#1}}
\nc\pro[1]{\begin{prop}
\label{#1}}
\nc\thm[1]{\begin{theorem}
\label{#1}}
\nc\teiri[1]{\begin{jtheorem}
\label{#1}}
\nc\cor[1]{\begin{corollary}
\label{#1}}
\nc\dfn[1]{\begin{defin}
\label{#1}}
\nc\sthm[1]{\begin{subth}
\label{#1}}
\nc\exm[1]{\begin{example}
\label{#1}
\begin{em}}
\nc\plm[1]{\begin{prblm}
\label{#1}
\begin{em}}
\nc\rmk[1]{\begin{remark}
\label{#1}
\begin{em}}
\nc\ntn[1]{\begin{notation}
\label{#1}
\begin{em}}
\nc\cau[1]{\begin{caution}
\label{#1}
\begin{em}}
\nc\imn[1]{\begin{importnota}
\label{#1}
\begin{em}}
\nc\cax[1]{\begin{cauex}
\label{#1}
\begin{em}}
\nc\con[1]{\begin{construction}
\label{#1}
\begin{em}}
\nc\ssthm[1]{\begin{ssubth}
\label{#1}
\begin{em}}
\nc\cnc[1]{\begin{conclusion}
\label{#1}
\begin{em}}

\nc\elem{\end{lemma}}
\nc\erest{\end{restroom}}
\nc\econj{\end{conjecture}}
\nc\eprodef{\end{propositiondef}}
\nc\eprt{\end{proto}}
\nc\ethm{\end{theorem}}
\nc\eteiri{\end{jtheorem}}
\nc\ecor{\end{corollary}}
\nc\edfn{\end{defin}}
\nc\esthm{\end{subth}}
\nc\epro{\end{prop}}
\nc\etri{\end{triviality}}
\nc\eexm{\end{em}
\end{example}}
\nc\ermk{\end{em}
\end{remark}}
\nc\eplm{\end{em}
\end{prblm}}
\nc\ecau{\end{em}
\end{caution}}
\nc\ecax{\end{em}
\end{cauex}}
\nc\eimn{\end{em}
\end{importnota}}
\nc\entn{\end{em}
\end{notation}}
\nc\econ{\end{em}
\end{construction}}
\nc\ecnc{\end{em}
\end{conclusion}}
\nc\essthm{\end{em}
\end{ssubth}}

\nc\GS{{\mathfrak S}}
\nc\pp{{\mathbb P}}
\nc\qq{{\mathbb Q}}
\nc\oo{{\mathcal O}}
\nc\sss{{\mathbb T}}
\nc\zz{{\mathbb Z}}
\nc\rr{{\mathbb R}}
\nc\cc{{\mathbb C}}
\nc\ff{{\mathbb F}}
\nc\ee{{\mathbb E}}
\nc\AAA{{\mathcal A}}
\nc\disjoint{{\amalg}}
\nc\id{{\rm id}}
\nc\vv{{\mathbb V}}
\nc\Ker{{\rm Ker }}
\nc\sgn{{\rm sgn}}
 \nc\CH{{\rm CH}}
\nc\Chow{{\mathcal C}}

\author{E. Javier Elizondo}
\address{Instituto de Matem\'aticas\\Universidad Nacional Aut\'onoma
  de M\'exico, Ciudad Universitaria, M\'exico
  DF 04510}
\email{javier@math.unam.mx}
\thanks{Partially support by program JSPS-CONACYT}
\author{Shun-ichi Kimura}
\address{Department of Mathematics, Graduate School of Science,
  Hiroshima University, Higashi-Hiroshima 739-8526, Japan}
\email{kimura@math.sci.hiroshima-u.ac.jp}

\title{Irrationality of Motivic series of Chow varieties}

\subjclass[2000]{14}
\keywords{Chow Motives, Chow varieties, Motivic Zeta}

\begin{document}
\maketitle
\section{Introduction}

The Euler characteristic of all the
Chow varieties, of a fixed projective variety,  can be collected  in a
formal power
series called the Euler-Chow series.  This series
coincides with the Hilbert series when the Picard group is a finite
generated free abelian group. It is an interesting open problem
to find for which varieties this series is rational. A few cases have
been computed, and it is suspected that the series is not rational for
the blow up of $\pp^2$ at nine points in general position. 

It is very natural to extend this series to Chow motives and ask the
question if the series is rational or to find a counterexample. In
this short paper we
generalized the series and show by  an example that the series
is not rational. This opens the question of what is
the geometrical meaning of the Euler-Chow series.

We would like to thank James D. Lewis for his suggestions in a
preliminary version.

 \section{Homological Chow Motives}

In this section we recall some definitions that we need in the next
sections.  

 \dfn{1-0} When $X$ and $Y$ are smooth complete varieties over $\cc$,
 an element $\alpha\in \CH_*(X\times Y)$ is called {\em a
 correspondence} from $X$ to $Y$, and written as $\alpha: X\vdash Y$.
 When $\alpha:X\vdash Y$ and $\beta: Y\vdash Z$ are correspondences,
 we define their composition $\beta\circ \alpha:X\vdash Z$ by
 $$\beta\circ\alpha=\pi_{XZ*}(\pi_{XY}^*\alpha\cdot \pi_{YZ}^*\beta)$$
 where $\pi$ are projections from $X\times Y\times Z$ to the products
 of indicated components, and $\underline{\ \ }\cdot \underline{\ \ }$
 denotes the intersection product.

Let $X=\disjoint X_i$ be the decomposition of $X$ into the irreducible
 (hence connected) components.  Then via the canonical isomorphism
 $\CH_*(X\times Y)\simeq \oplus \CH_*(X_i\times Y)$, a correspondence
 $\alpha: X\vdash Y$ decomposes into $\alpha=\oplus \alpha_i$, with
 $\alpha_i\in \CH_*(X_i\times Y)$.  A correspondence $\alpha: X\vdash
 Y$ is said to have {\em relative dimension $d$} when $\alpha_i\in
 \CH_{d+\dim X_i}(X_i\times Y)$. \edfn 

 \dfn{1-1} A {\em  (homological)  Chow motive} is a triple $(X, p, n)$
 where $X$ is a smooth complete scheme over $\cc$,\,  $p: X\vdash X$
 is a correspondence with relative dimension $0$ such that $p\circ p=p$
 (namely $p$ is idempotent), and $n$ is an integer.  The Chow motive $(X,
 [\Delta_X], 0)$ is called the {\em \bf Chow motive of $\mathbf X$}
 and it is written as $h(X)$, where $\Delta_X\subset X\times X$ is the
 diagonal subvariety. 

 When $M=(X, p, n)$ and $N=(Y, q, m)$ are Chow motives, then a {\em \bf
 morphism of Chow motives} from $M$ to $N$ is a correspondence
 $\alpha: X\vdash Y$ such that $\alpha$ has a relative dimension $n-m$
 and $q\circ \alpha\circ p=\alpha$. Finally, the category of Chow motives is
 denoted as ${\mathcal ChMot}$. 

 When $f: X\to Y$ is a morphism of smooth complete schemes, then its
 graph $[\Gamma_f]\in\CH_*(X\times Y)$ determines a morphism $h(X)\to
 h(Y)$, and this morphism is written as $h(f)$.

 For a motive $M=(X, p, n)$, we define its $i$-th homology $H_i(M)$ to
 be the image of the projector 
$$H_i(M):=p_*(H_{i+2n}(X, \qq))\subset H_{i+2n}(X, \qq).$$  

 {\em \bf The tensor product} $M \otimes N$ of two Chow motives $M=(X,
 p,n)$ and $N=(Y, q, m)$ is defined to be 
$$M\otimes N:=(X\times Y, \pi_X^*p\cdot \pi_Y^*q, n+m)$$ 
where $\pi_X: (X\times Y)\times (X\times Y)\to X\times X$ and $\pi_Y:
(X\times Y)\times (X\times Y)\to Y\times Y$ are projections, and
``$\cdot$'' is the intersection product.  
\edfn 

 \rmk{1-2} $h$ determines a covariant functor from the category of
 smooth complete schemes  to the category of Chow motives. \ermk

 \rmk{1-3} Classically, Chow motives are defined cohomologically, and
 the cohomology functor $h$ is a contravariant functor
 (cf. \cite{scholl-motives}).  In this paper, we need homological
 operation on the 
 Chow motives rather than cohomological one, so we choose to use
 homological (and hence covariant) definition of Chow motives.  \ermk

\lem{1-4} Let $P\in \pp^1$ be some point in $\pp^1$, and consider a
point $Pt\,: =\,
  {\rm Spec} \,(\cc)$ outside of $\pp^1$. Then, in the
  category of Chow motives, we have an isomorphism 
$$(\pp^1, [P\times \pp^1], -1)\simeq (Pt, [\Delta_{Pt}], 0).$$ 
Furthermore, for any motive $M$ there is an isomorphism  
$$M\otimes (Pt, [\Delta_{Pt}], 0)\simeq M.$$  \elem 

\begin{proof} The morphisms $[P\times Pt]: \pp^1\vdash Pt$ and
  $[Pt\times \pp^1]: Pt\vdash \pp^1$ give the isomorphism between
  $(\pp^1, [P\times \pp^1], -1)$ and $(Pt, [\Delta_{Pt}], 0)$.  The
  isomorphism $M\otimes (Pt, [\Delta_{Pt}], 0)\simeq M$ follows
  immediately from the definition of tensor product.\end{proof}

\cor{1-5} We have an isomorphism of Chow motives 
$$\left( X, p, n\right)\simeq \left( X\times (\pp^1)^{n-m}, p\times [(P\times
\pp^1)^{n-m}], m\right), $$ 
where $P\in \pp^1$ is any fixed point. \ecor 

\dfn{1-6}  Let $M=(X, p, n)$ and $N=(Y, q, m)$ be Chow motives.  When
$n=m$, we define $M\oplus N$ by 
$$M\oplus N :=(X\disjoint Y, p\, \disjoint \, q,\, n).$$
When $n\not=m$, say $n>m$, we define $M\oplus N$ by 
$$M\oplus N:=(X\times (\pp^1)^{n-m}\disjoint Y,  p\times [(P\times
\pp^1)^{n-m}]\disjoint q, m).$$ \edfn

\dfn{1-7} We define \, $K_0({\mathcal ChMot})$ \,to be \,$({\mathcal
  ChMot}\times {\mathcal ChMot})/\sim$\,  where \,$(M_1, N_1)\sim (M_2,
  N_2)$\,  iff \, $M_1\oplus N_2\simeq M_2\oplus N_1$.
We write the object $(M, N)$ as $[M]-[N]$. \edfn

\rmk{1-8} Usual argument says that $K_0({\mathcal ChMot})$ is a ring
with $\oplus$ as addition and $\otimes$ as multiplication.  \ermk

\rmk{1-9} The cohomology series extends to $K_0({\mathcal ChMot})$ by
sending $[M]-[N]$ to $\sum (\dim H_i(M)-\dim H_i(N))t^i$.  This gives
a ring homomorphism from $K_0({\mathcal ChMot})\to \zz[t]$.
By further substituting $t=-1$, there is a ring homomorphism
$K_0({\mathcal ChMot})\to \zz$, which sends $h(X)$ to $\chi(X)$, the Euler characteristic of $X$.
\ermk

\section{Chow Series of Chow Varieties}

In this section we define the {\it Chow series} for Chow motives
which generalizes the {\it Euler-Chow series}  for
Chow varieties defined in \cite{eli-tor} and \cite{eli-pau}. 

\dfn{2-1} Let $X$ be a smooth projective scheme, and $\lambda\in
H_{2p}(X, \zz)$ be its homology class.  Define ${\Chow_\lambda}(X)$ to
be the space of all effective cycles on $X$ whose homology class is $\lambda$.
We define $\Chow_p(X)$ to be the disjoint union of all
$\Chow_\lambda(X)$ with $\lambda\in H_{2p}(X, \zz)$. \edfn

\rmk{2-2} The space $\Chow_\lambda(X)$ is an open and closed subscheme
of the classical Chow variety $\Chow_{p, d}(X)$ where $d$ is the
degree of $\lambda$ for the given projective structure, hence
$\Chow_{\lambda}(X)$ is a projective scheme.  See \cite[Lemma 1.1]{eli-tor}.\ermk

\rmk{2-3} The space of all effective $p$-cycles $\Chow_p(X)$ has the
canonical monoid structure, $(\alpha, \beta)\mapsto
\alpha+\beta$, the addition of $p$-cycles, which is compatible with
the addition of 
the homology class: $\Chow_\lambda(X)\times \Chow_\mu(X)\to
\Chow_{\lambda+\mu}(X)$. \ermk

\dfn{2-4} Let $C$ be the monoid of homology classes of effective
$p$-cycles in $H_{2p}(X, \zz)$, and $K_0({\mathcal ChMod})[[C]]$ be
the set of all $K_0({\mathcal ChMod})$ valued functions on $C$. \edfn

\rmk{2-5} $K_0({\mathcal ChMot})[[C]]$ has a ring structure with
multiplication defined by the convolution product, i.e., for $f, g\in
K_0({\mathcal ChMot})[[C]]$, we define $$(f\times
g)(\lambda):=\sum_{\lambda=\mu_1+\mu_2}f(\mu_1)\otimes f(\mu_2)$$ 
See \cite[Prop. 1.3]{eli-tor} for the well-definedness of this product. \ermk

\dfn{2-6} Let $X$ be a projective variety.  We define the {\em Chow
  series} of $X$, in dimension $p$, to be the element $C_p(X)\in
K_0({\mathcal ChMot})[[C]]$ by $C_p(X)(\lambda)=[\Chow_\lambda(X)]$
for $\lambda\in H_{2p}(X, \zz)$.  By convention, if $\Chow_\lambda$ is
the empty set $\phi$, then we define $[\phi]=0$ in $K_0({\mathcal
  ChMot})$. 
\edfn

\dfn{2-7} Let $K_0({\mathcal ChMot})[C]$ be the monoid-ring of $C$
over $K_0({\mathcal ChMot})$.  This ring consists of all elements of
$K_0({\mathcal ChMot})[[C]]$ with finite support, and sometimes we
referto an element of  $K_0({\mathcal ChMot})[C]$ as a polynomial.  \edfn

\dfn{2-8} An element  $\varphi\in K_0({\mathcal ChMot})[[C]]$ is {\em
  rational} if for some polynomials $f, g\in K_0({\mathcal ChMot})[C]$
with $f$ invertible in $K_0({\mathcal ChMot})[[C]]\otimes \qq$, we
have $f\varphi=g$ in $K_0({\mathcal ChMot})[[C]]$.   \edfn

\rmk{2-8-1} It is not known if $K_0({\mathcal ChMot})$ is an integral
domain.  A polynomial $f\in K_0({\mathcal ChMot})[C]$ is invertible in
$K_0({\mathcal ChMot})[[C]]\otimes \qq$ iff $f(0)$ is invertible in
$K_0({\mathcal ChMot})\otimes \qq$, where $0\in H_{2p}(X, \zz)$ is the
zero object.  \ermk

\rmk{2-9} If we assume the finite dimensional conjecture (see
\cite{kimura-findim}), 
or the stronger assertion of the Bloch-Beilinson
conjecture (see \cite{jannsen-filtration}), then for $p=0$, the Chow
series $\Chow_0(X)$ is 
rational (see \cite{andre-motidimension}).  Also by using the ring
homomorphism from 
$K_0({\mathcal ChMot})\to \zz$ which sends $h(X)$ to $\chi(X)$ (see
Remark \ref{1-9}), if we know that $\Chow_p(X)$ is rational, then it
implies that the Euler-chow Series $E_p(X)=\sum_\lambda
\chi(\Chow_\lambda)\lambda$ is also rational. The series is known to
be rational for toric varieties (see \cite{eli-tor}), particular cases
for Grassmannians (see \cite{eli-pau}) and some Del Pezzo surfaces.  \ermk 

\rmk{2-10} The rational elements of $K_0({\mathcal ChMot})[[C]]$ is
closed under addition and multiplication.  Also if a rational element
is invertible, then its inverse is also rational.  \ermk

\newpage

\section{Chow series of projective spaces}

In this section, we prove that the Chow series $C_{n-1}(\pp^n)$ is not
rational, if $n>1$. 

\lem{3-1} Let $CS: K_0({\mathcal ChMot})[[C]]\to \zz[t][[C]]$ be the
cohomology series homomorphism. 
If $CS(\varphi)$ is not a rational function in $\qq[t][[C]]$ for
$\varphi\in K_0({\mathcal ChMot})[[C]]$, then $\varphi$ is not
rational.  \elem

\begin{proof} If $\varphi$ is rational, then there are polynomials $f,
  g\in K_0({\mathcal ChMot})[C]$ with $f$ invertible in $K_0({\mathcal
    ChMot})[[C]]\otimes \qq$ such that $f\varphi=g$.  We send both
  sides by the ring homomorphism $CS$ to find that
  $CS(f)CS(\varphi)=CS(g)$, where $CS(f)$ and $CS(g)$ are polynomials
  in $\zz[t][C]$.  Let $h\in  K_0({\mathcal ChMot})[[C]]\otimes \qq$
  be the inverse of $f$, then $CS(f)CS(g)=1$ in $\qq[t][[C]]$.  In
  particular, the constant term of $CS(f)$ is not zero, and hence non
  zero divisor, and we can write $CS(\varphi)=CS(g)/CS(f)$ as a
  rational function, a contradiction.  \end{proof} 

\dfn{3-2} Let $\varphi\in \zz[[s, t]]$ be a power series.  We write
$\varphi$ as $$\varphi(s, t)=\sum_{i=0}^\infty a_i(t)s^i$$ where
$a_i(t)\in \zz[[t]]$ for each $i$.  We say that $\varphi$ has a gap
sequence $\{\alpha_1, \alpha_2, \alpha_3, \ldots\}$, with $a_i \in \N$\, 
$\forall i$, if there is a sequence of natural numbers
$d_1<d_2<d_3<\cdots$ such that $a_{d_i}(t)\not=0$ and
$a_{d_i+1}(t)=a_{d_i+2}(t)=\cdots=a_{d_i+\alpha_i}(t)=0$.  \edfn

\lem{3-3} If $\varphi\in  \zz[[s, t]]$ has a gap sequence $\{\alpha_1,
\alpha_2, \alpha_3, \ldots\}$ such  that $\displaystyle \lim_{n\to
  \infty}\alpha_n=\infty$, then $\varphi$ is not a rational function.
Namely, $\varphi$ is not in $\qq(s, t)$.   \elem

\begin{proof} Assume that $\varphi(s, t)f(s, t)=g(s, t)$ for some
  polynomials $f(s,  t)$ and $g(s, t)$ in $\zz[s, t]$.  By the
  definition of the gap sequence, we can take $d_1<d_2<d_3<\cdots$.
  When $i$ is large enough, we have $\alpha_i>\deg_sf$ and
  $d_i>\deg_sg$.  Then term of $\varphi(s, t)f(s, t)$ with the degree
  in $s$ to be $\deg_sf+d_i$ is non zero, and on  the other hand, the
  term of $g(s, t)$ with the degree in $s$ to be $\deg_sf+d_i$ is
  zero, a contradiction.   \end{proof}

\pro{3-4} The cohomology series   of $\Chow_{n-1}(\pp^n)$  is given by
$$
\sum_i \dim
H_i(\Chow_{n-1}(\pp^n))t^i\, = \, 
\sum\frac{1-s^{2\left(\begin{smallmatrix}d+n\\d\end{smallmatrix}\right)}}{1-s^2}\,\,t^d
$$.
\epro

\begin{proof} $H_{2n-2}(\pp^n, \zz)\simeq \zz$, and one can identify
  $C=\{0, 1, 2, \ldots\}$, by degree.  The degree $d$ part of the Chow
  variety is $\Chow_{n-1, d}(\pp^n)\simeq
  \pp^{\left(\begin{smallmatrix}d+n\\d\end{smallmatrix}\right)-1}$, from which
  the proposition immediately follows.  \end{proof}

\thm{3-5} The Chow series $\sum [\Chow_{n-1, d}]t^d\in K_0({\mathcal
  ChMot})[[C]]$ is not rational.  \ethm

\begin{proof} By Lemma \ref{3-1}, it is enough to show that
$$
\varphi:=
\sum\frac{1-s^{2\left(\begin{smallmatrix}d+n\\d\end{smallmatrix}\right)}}{1-s^2}t^d 
$$  
  is not rational in $\qq(s, t)$.  Assume that $\varphi$ is a rational
  function, then
  $$
(1-s^2)\varphi=\sum\left(1-s^{2\left(\begin{smallmatrix}d+n\\d\end{smallmatrix}\right)}\right)t^d
$$
  is also a rational function.  But by setting
  $d_i:=2\left(\begin{smallmatrix}d+n\\d\end{smallmatrix}\right)$,
  $(1-s^2)\varphi$ has the gap sequence $\{\alpha_1, \alpha_2,
  \ldots\}$ such that
 $$\alpha_d=2\left(\begin{array}{cc}d+1+n\\d+1\end{array}\right)-2
  \left(\begin{array}{cc}d+n\\d
\end{array}\right)=2\left(\begin{array}{cc}d+n\\n-1\end{array}\right)
$$
  which tends to infinity.  By Lemma \ref{3-3}, it implies that
  $\varphi$ is not rational.   

\end{proof}

\end{document}